\documentclass[11pt,reqno]{amsart}
\usepackage[letterpaper,margin=1.15in]{geometry}
\usepackage{amsmath,amssymb,mathtools}
\usepackage{graphicx,xcolor}
\usepackage{hyperref}
\numberwithin{equation}{section}
\newtheorem{theorem}{Theorem}[section]
\theoremstyle{definition}
\newtheorem{thmA}{Theorem}

\newtheorem{thmB}{Theorem}

\newtheorem*{problemA}{Problem A}
\newtheorem*{problemB}{Problem B}

\theoremstyle{plain}

\newtheorem{open problem}{Open Problem}[section]
\newtheorem{open question}{Open Quesion}[section]
\newtheorem{proposition}[theorem]{Proposition}
\newtheorem{lemma}[theorem]{Lemma}
\newtheorem{corollary}[theorem]{Corollary}
\theoremstyle{definition}
\newtheorem{remark}[theorem]{Remark}
\newcommand{\R}{\mathbb R}
\allowdisplaybreaks[2]
\title[How to arrange four points in the plane]{How to arrange four points in the plane}
\author{Senping Luo}
\address{School of Mathematics and Statistics, Jiangxi Normal University, Nanchang 330022, China}
\email{luosp1989@163.com}
\begin{document}
\begin{abstract}
Arranging points in the plane appears to be an elementary geometric problem,
yet it is deeply rooted in statistical physics. As a toy model, we consider the configuration
of four points. We completely classify the minimal configurations of four points interaction
energy in $\R^2$ under the general Lennard-Jones potentials $V(r)=r^{-p}-br^{-q}$, where $p>q>0$
and $b\in\R$. We find that there is a threshold curve
\[
 D_{p,q}=2^{-(p-q)/2}\Bigl(1+\frac{p-q}{q+2}\bigl(1+2^{-(q+2)/2}\bigr)\Bigr)
\]
to determine the minimal shapes. The main results can be stated concisely as follows: consider
$\min\sum_{1\le i<j\le4}\bigl(|x_i-x_j|^{-p}-b|x_i-x_j|^{-q}\bigr)$.
For $b\le0$, the minimizer does not exist. For $b>0$, the minimizer is
unique up to rotation and translation. Furthermore, the minimizer is
\[
 \begin{cases}
 \text{a square},&\hbox{if}\;\; D_{p,q}\ge1,\\
 \text{a rhombus with angle }\theta\in(\pi/3,\pi/2),&\hbox{if}\;\; D_{p,q}<1.
 \end{cases}
\]
\end{abstract}
\maketitle

\section{Introduction and statement of main results}\label{Section1}

At low temperatures, atoms arrange themselves into periodic structures and regular shapes.
How interatomic forces the periodic structures and regular shapes is a fundamental problem
in statistical and solid mechanics~\cite{DLF,Radin1987,BL}. To model the physical problems from the first principle and simplify the statements, we formulate the problems in a purely mathematical form, and we call that how to arrange points minimizing certain functionals.

Let $X=(x_1,\ldots,x_N)\in(\R^2)^N$ consist of $N\ge1$ distinct
points. For a radial two-body potential $V$, the $N$ points interaction energy is
\begin{equation}\label{1Equation1}
 E_V^{(N)}(X)=\sum_{1\le i<j\le N} V\bigl(|x_i-x_j|\bigr)
 =\frac12\sum_{\substack{1\le i,j\le N\\ i\ne j}}V\bigl(|x_i-x_j|\bigr).
\end{equation}
The $N$ points interaction energy arises in lots of physical and mathematical problems, for which we refer to the recent Lecture notes by Serfaty \cite{Serfaty}, and the reviews by Radin \cite{Radin1987} and Blanc-Lewin \cite{BL}.

There are freedom on particle numbers $N$ and the potentials $V$. We start from a simple yet useful potential.
The Heitmann--Radin potential is defined as follows:
\begin{equation}\label{HRP}
 V_{\rm HR}(r)=
 \begin{cases}
  +\infty,&0<r<1,\\
  -1,&r=1,\\
  0,&r>1.
 \end{cases}
\end{equation}
Correspondingly, the Heitmann--Radin energy~\cite{HR} is $N$ points interaction energy \eqref{1Equation1}
$
 E_{V_{\rm HR}}^{(N)}(X)
$
under the potential \eqref{HRP}.

In 1980, Heitmann--Radin proves the following celebrated result.
\begin{thmA}[Heitmann--Radin \cite{HR}]\label{TheoremA}
Every minimizer $X$ of $E_{V_{\rm HR}}^{(N)}(X)$ is subset of triangular lattice. That is,
\[
 \{x_1,\ldots,x_N\}
 \subset y+Q\left(\mathbb Z(1,0)+\mathbb Z\bigl(\tfrac12,\tfrac{\sqrt3}2\bigr)\right)
\]
for some $y\in\R^2$ and $Q\in\mathrm{SO}(2)$.
\end{thmA}
A new proof of Theorem~\ref{TheoremA} can be found in De Luca-Friesecke \cite{DFGaussBonnet}. Theorem~\ref{TheoremA} does not include any uniqueness of the minimizers
up to translation and rotation. This question was completely solved by De Luca-Friesecke \cite{DLF} as follows:
\begin{thmB}[De Luca-Friesecke \cite{DLF}]\label{TheoremB}
Let $N\in\mathbb{N}$. The minimizers of $E_{V_{\rm HR}}^{(N)}(X)$  are unique up to translation and rotation
if and only if either
\begin{enumerate}
\renewcommand{\labelenumi}{\textup{(\alph{enumi})}}
\item $N=3s^2+3s+1$, for some $s\in\mathbb N\cup\{0\}$
\item $N=3s^2+3s+1+(s+1)k+s$, for some $s\in\mathbb N\cup\{0\}$ and some $k\in\{0,1,2,3,4\}$.
\end{enumerate}
\end{thmB}
Theorems~\ref{TheoremA} and~\ref{TheoremB} together provide a complete classification including existence and uniqueness on the minimizers of the $N$ points interaction energy under Heitmann--Radin potential.

The Heitmann--Radin potential \eqref{HRP} is useful to model solid collision among particles, while physically there is a widely used Lennard-Jones potential as follows:

\begin{equation}\label{VLJ}
V_{\rm LJ}(r):=r^{-2p}-2r^{-p},\;\; p>0.
\end{equation}
See, for example, in \cite{DFGaussBonnet}. The Lennard-Jones has long history and receives extensive studies (\cite{Mie,Jones1924,LSH2024,LWJEMS,BL}).
As noted in De Luca-Friesecke \cite{DLF}, the Lennard-Jones potential converges to Heitmann--Radin potential as $p\rightarrow+\infty$.
This reveals a connection between these two potentials.
To fully understand particles interaction under Lennard-Jones potentials, it is natural to ask the following problem, in parallel to Theorems~\ref{TheoremA} and~\ref{TheoremB}:
\begin{problemA}
Classify the minimizers of $N$-points interaction energy \eqref{1Equation1} under the Lennard-Jones potential \eqref{VLJ}.
\end{problemA}

There is a related interesting problem concerning Theorem~\ref{TheoremA}. The triangular and hexagonal shapes are prevail and ubiquitous in nature. Understanding the triangular and hexagonal shapes from mathematical perspective from the first principle is important for both in physics and applied and pure mathematics. There are many interesting works; see Montgomery \cite{Montgomery} for the two-dimensional theta functions,
 B\'eterin-Faulhuber-Steinerberger \cite{BFS} for shifted Gaussian lattice sum with wide applications, Chen-Oshita \cite{CO} for copolymers, Sandier-Serfaty \cite{SS2012} for Ginzburg-Landau superconductivity, Theil \cite{Theil} and E-Li \cite{E} for triangular crystallization and hexagonal crystallization, respectively.

 While, we shall raise the following problem from the opposite direction.
\begin{problemB}
For what situations,
the minimizers of $N$-points interaction energy \eqref{1Equation1} are not the triangular or regular shapes (or their subsets).
\end{problemB}

For $N=4$, Theorem~\ref{TheoremA} indicates that the rhombus with angle $\pi/3$ (two equilateral triangles) minimizes the four points interaction energy under Heitmann--Radin potential.

To answer Problems A and B, we focus on $N=4$ as a toy model and consider the most general form of Lennard-Jones potentials $V(r)=r^{-p}-br^{-q}$ with $p>q>0$ and $b\in\R$. Then the four points interaction energy has the form as follows:
\begin{equation}\label{1Equation3}
 E_{p,q,b}(X):
 =\sum_{1\le i<j\le4}\bigl(|x_i-x_j|^{-p}-b|x_i-x_j|^{-q}\bigr).
\end{equation}

We provide a complete classification on the minimizers of this four points interaction energy for all
$p>q>0$ and $b\in\R$.
It is stated as follows:
\begin{theorem}[Minimal configurations]\label{1Theorem1}
Let $p>q>0$ and $b\in\R$. Consider the four points interaction energy under the general Lennard-Jones potentials as follows:
\begin{equation}\nonumber
\sum_{1\le i<j\le4}\bigl(|x_i-x_j|^{-p}-b|x_i-x_j|^{-q}\bigr).
\end{equation}
Then, there is a threshold curve defined as follows:
\begin{equation}\label{1Equation4}
 D_{p,q}
 :=2^{-(p-q)/2}\Bigl(1+\frac{p-q}{q+2}
                            \bigl(1+2^{-(q+2)/2}\bigr)\Bigr).
\end{equation}
The following hold.

\begin{enumerate}
\renewcommand{\labelenumi}{\textup{(\roman{enumi})}}
\item For $b\le0$, the minimizer does not exist.

\item For $b>0$, the minimizer is a rhombus.
Up to rotation and translation, it is
\[
 X_{p,q,b}^*
 =a_b\begin{pmatrix}
 \cos(\theta_{p,q}/2)&0&-\cos(\theta_{p,q}/2)&0\\
 0&\sin(\theta_{p,q}/2)&0&-\sin(\theta_{p,q}/2)
 \end{pmatrix},
\qquad a_b>0.
\]

\item For $b>0$, the angle $ \theta_{p,q}$ of the rhombus satisfies
\[
 \begin{cases}
 \theta_{p,q}=\pi/2,& D_{p,q}\ge1,\\
 \pi/3<\theta_{p,q}<\pi/2,& D_{p,q}<1 .
 \end{cases}
\]
In particular, the minimizers are rhombi whose angles can never be $\pi/3$.

\item For $b>0$, the following hold.
\begin{enumerate}
\renewcommand{\labelenumii}{\textup{(\alph{enumii})}}
\item The map $(p,q)\mapsto\theta_{p,q}$ is real analytic on
$\{(p,q):p>q>0,\ D_{p,q}<1\}$. Here,
$\partial_p\theta_{p,q},\partial_q\theta_{p,q}<0$.

\item For each fixed $q>0$,
\[
 \lim_{p\to\infty}(p+2)\bigl(\theta_{p,q}-\pi/3\bigr)
 =\frac2{\sqrt3}\log\frac{1+3^{-(q+2)/2}}{1-3^{-(q+2)/2}}.
\]

\item For $p=2q$. As $q\to\infty$,
\[
 \theta_{2q,q}=\frac\pi3+\frac{4\cdot3^{-q/2}}{3\sqrt3\,q}
                         +O(3^{-q}/q).
\]
\item If $p^*>q^*>0$ and $D_{p^*,q^*}=1$, then
\[
 \lim_{\substack{(p,q)\to(p^*,q^*)\\D_{p,q}<1}}
 \frac{\pi/2-\theta_{p,q}}{\sqrt{1-D_{p,q}}}
 =
 \frac{2\sqrt6}
 {\sqrt{(p^*-q^*)\bigl[p^*+q^*+10+(4-2p^*+q^*)2^{-(p^*+2)/2}\bigr]}}.
\]
\end{enumerate}
\end{enumerate}
\end{theorem}

\begin{remark}\label{1Remark2}
For $p=2q$ and $b=2$, the minimizer in \textup{(ii)} converges to a
rhombus with angle $\pi/3$. More precisely,
\[
 X_{2q,q,2}^*\longrightarrow
 \begin{pmatrix}
  \sqrt3/2&0&-\sqrt3/2&0\\
  0&1/2&0&-1/2
 \end{pmatrix}
 \qquad (q\to\infty).
\]
The limit gives the minimizer of $E_{V_{\rm HR}}^{(4)}(X)$.
\end{remark}

 The square shape is not the first time being proved as the minimizer; see
B\'eterin-De Luca-Petrache \cite{BDP} on square crystallization.
Recently, Del Nin-De Luca proved square crystallization
for supremum-norm sticky disks~\cite{DNLDSquare}. The square lattice shape appear also in triblock copolymer (Luo-Ren-Wei \cite{LRW}) and Bose-Einstein condensates (Mueller-Ho \cite{Mue2002} and Luo-Wei \cite{LWMH}).

\begin{figure}[htbp]
\centering
\includegraphics[width=0.5\linewidth]{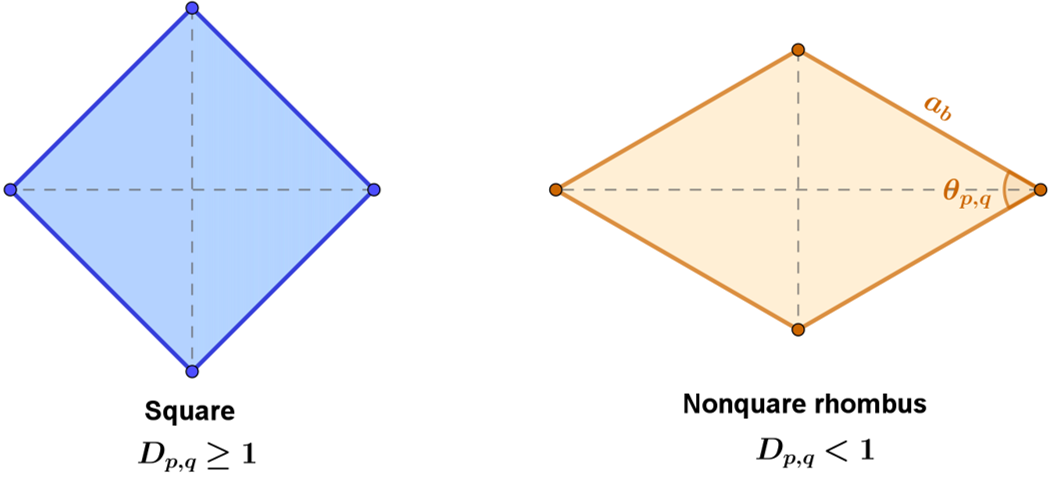}
\caption{The minimal shapes.}\label{1Figure1}
\end{figure}

Considering thermodynamic of interaction energy \eqref{1Equation1} per particle, and assuming the particles arrange themselves into a lattice, we obtain the $f-$potential lattice energy
 \begin{equation}\aligned\label{EFL}
E_f(\Lambda):=\sum_{P\in \Lambda \backslash\{0\}} f\bigl(|P|^2\bigr),\; |\cdot|\;\hbox{is the Euclidean norm on}\;\mathbb{R}^2.
\endaligned\end{equation}
\begin{figure}[htbp]
\centering
\includegraphics[width=0.5\linewidth]{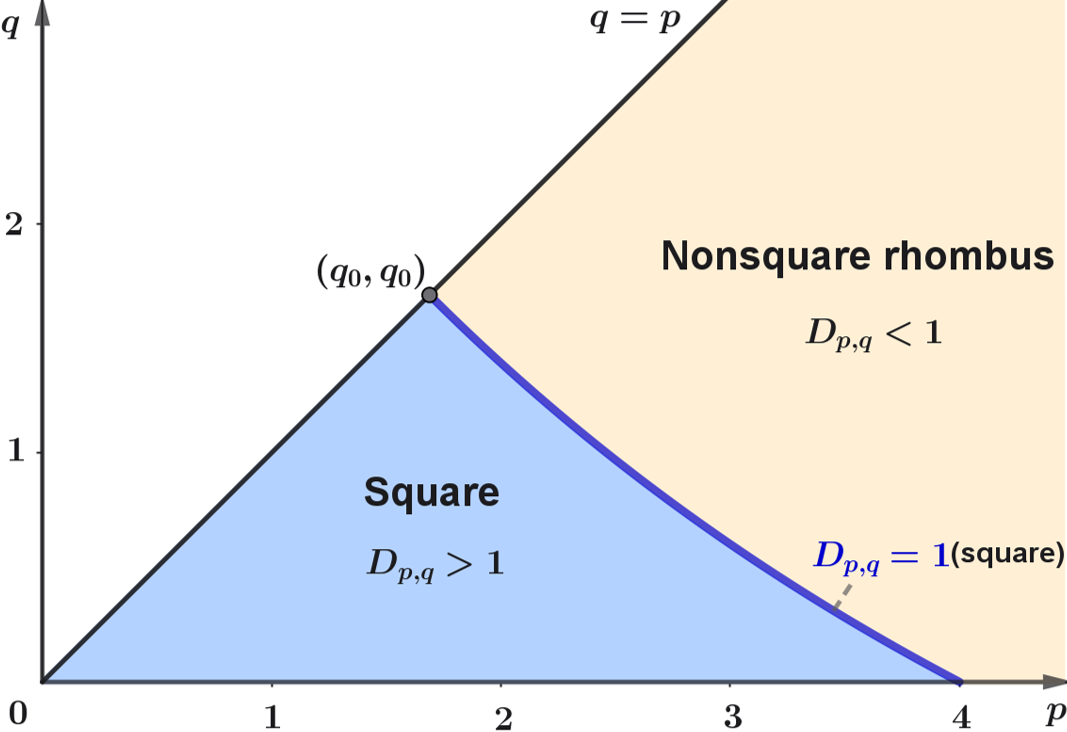}
\caption{The phase diagram for $p>q>0$.
The endpoints of $D_{p,q}=1$ are $(4,0)$ and $(q_0,q_0)$, where $q_0\approx1.68887$.}\label{1Figure2}
\end{figure}

Here $\Lambda$ denotes a lattice. This is another related interesting topic on lattices. We refer to B\'eterin-Petrache \cite{BP} for a detailed review on lattices and Blanc-Lewin \cite{BL} for the broad review on crystallization including lattice minimizations and special functions as a section.
Cohn-Kumar-Miller-Radchenko-Viazovska \cite{CKMRV} prove the universality of $E_8$ and Leech lattices in dimensions 8 and 24, respectively.
Interestingly, Petrache-Serfaty \cite{PS2020} then prove the crystallization results in these two dimensions building on this remarkable theorem (\cite{CKMRV}). For the $f-$potential lattice energy \eqref{EFL} under the classical Lennard-Jones potential in dimension two, we refer to B\'eterin and his collaborators \cite{Bet2018,BZ,Betermin2021a} for pioneering work and Luo-Wei \cite{LWJEMS} for a complete resolution; while this problem is open in dimensions higher than two.

The paper is organized as follows:
Section~\ref{Section2} reduces the minimizers to rhombi. It rules out
nonconvex minimizers by comparison with squares. For a convex
minimizer, Connelly's rank-one equilibrium stress~\cite{Connelly}
shows that the diagonals are perpendicular and bisect each other.
We consider the minimization over rhombi in Section~\ref{Section3}.
We consider the classical Lennard-Jones case in Section~\ref{Section4}.

\section{Reduction to rhombi}\label{Section2}

For $b\leq0$, dilation gives a zero infimum that is not attained.
For $b>0$, \eqref{2Equation3} and Lemmas~\ref{2Lemma1}--\ref{2Lemma3}
give existence and reduce the problem to rhombi.

Write
\[
 e_{p,q,b}:=\inf_X E_{p,q,b}(X),\qquad r_{ij}:=|x_i-x_j|.
\]
For $b\le0$,
\[
 0<E_{p,q,b}(X),\qquad \lim_{t\to\infty}E_{p,q,b}(tX)=0,\qquad
 \Bigl(\min_{1\le i<j\le4}r_{ij}\Bigr)^{-p}\le E_{p,q,b}(X).
\]
The energy tends to zero exactly when $\min_{1\le i<j\le4}r_{ij}\to\infty$.
Differentiation under dilation gives
\[
 \left.\frac{d}{dt}E_{p,q,b}(tX)\right|_{t=1}
 =-p\sum_{1\le i<j\le4}|x_i-x_j|^{-p}
   +bq\sum_{1\le i<j\le4}|x_i-x_j|^{-q}<0,
\]
so there are no critical configurations. This proves part \textup{(i)} of Theorem~\ref{1Theorem1}.

For $b>0$, the potential $r^{-p}-br^{-q}$ has a unique minimum at $\ell_b$,
with well depth $\varepsilon_b$, where
\begin{equation}\label{2Equation1}
 \ell_b=\Bigl(\frac{p}{bq}\Bigr)^{1/(p-q)},\qquad
 \varepsilon_b=-\bigl(\ell_b^{-p}-b\ell_b^{-q}\bigr)
 =\frac{p-q}{q}\Bigl(\frac{bq}{p}\Bigr)^{p/(p-q)}.
\end{equation}
Define the normalized potential and energy by
\begin{equation}\label{2Equation2}
 W(r)=\frac{q r^{-p}-p r^{-q}}{p-q},\qquad
 \mathcal E_{p,q}(X)=\sum_{1\le i<j\le4}W(r_{ij}),
\end{equation}
so that
\begin{equation}\label{2Equation3}
 E_{p,q,b}(\ell_b X)=\varepsilon_b\mathcal E_{p,q}(X).
\end{equation}
Differentiating \eqref{2Equation2} gives
\begin{equation}\label{2Equation4}
 W'(r)=\frac{pq}{p-q}r^{-p-1}\bigl(r^{p-q}-1\bigr).
\end{equation}
Then
\begin{equation}\label{2Equation4a}
 W(r)\ge W(1)=-1.
\end{equation}
$W(r)<0$ exactly when $r>(q/p)^{1/(p-q)}$.

\begin{lemma}[Existence and compactness]\label{2Lemma1}
Let $p>q>0$ and $m_{p,q}=\inf_X\mathcal E_{p,q}(X)$. A
minimizer exists, and
\begin{equation}\label{2Equation5}
 m_{p,q}<-5+W(\sqrt3)<-5.
\end{equation}
Every distance in a minimizer lies in $\bigl((q/p)^{1/(p-q)},\sqrt3\bigr)$.
\end{lemma}
\begin{proof}
A rhombus of side length $a$ and acute angle $\pi/3$ has diagonals $a$ and
$\sqrt3\,a$, and hence energy
\[
 5W(a)+W(\sqrt3\,a).
\]
By \eqref{2Equation4},
\[
 \left.\frac{d}{da}\Bigl[5W(a)+W(\sqrt3\,a)\Bigr]\right|_{a=1}
 =\sqrt3\,W'(\sqrt3)>0.
\]
For $a<1$ sufficiently close to $1$, the energy is less than $-5+W(\sqrt3)$, its value at
$a=1$. Since \eqref{2Equation2} gives $W(\sqrt3)<0$, this proves \eqref{2Equation5}.

Fix $1\le i<j\le4$. Each of the other five pair terms is at least $-1$ by \eqref{2Equation4a}, so
\[
 \mathcal E_{p,q}(X)=W(r_{ij})+
 \sum_{\substack{1\le k<\ell\le4\\(k,\ell)\ne(i,j)}}W(r_{k\ell})
 \ge W(r_{ij})-5.
\]
If $\mathcal E_{p,q}(X)<-5+W(\sqrt3)$, the preceding estimate gives
\[
 W(r_{ij})<W(\sqrt3)<0\qquad (1\le i<j\le4).
\]
By \eqref{2Equation2}, $W(r_{ij})<0$ implies $r_{ij}>(q/p)^{1/(p-q)}$.
Since \eqref{2Equation4} makes $W$ strictly increasing on $[1,\infty)$,
$W(r_{ij})<W(\sqrt3)$ also implies $r_{ij}<\sqrt3$.

By \eqref{2Equation4a}, $m_{p,q}\ge-6$. Translation invariance of $\mathcal E_{p,q}$ allows us to impose $x_1^{(n)}=0$
on a minimizing sequence $X^{(n)}$. By \eqref{2Equation5} and the preceding distance bounds,
$|x_i^{(n)}|<\sqrt3$ for $2\le i\le4$ and all sufficiently large $n$.
A subsequence therefore converges to some $X^*=(x_1^*,\ldots,x_4^*)$.
Passing to the limit in the lower bound gives
\begin{equation}\label{21a}
 |x_i^*-x_j^*|\ge(q/p)^{1/(p-q)}>0\qquad(1\le i<j\le4),
\end{equation}
so continuity gives $\mathcal E_{p,q}(X^*)=m_{p,q}$. By \eqref{2Equation5}, every minimizer lies in
the strict sublevel considered above and hence satisfies both strict distance bounds.
\end{proof}

\begin{lemma}[Comparison with a square]\label{2Lemma2}
Let $p>q>0$. Suppose one point of $X$ lies in the
convex hull of the other three. Then there is a square $Y$ such that $\mathcal E_{p,q}(Y)<\mathcal E_{p,q}(X)$.
\end{lemma}
\begin{proof}
Set
\[
 F(t)=1+W(e^t),\qquad F_\lambda(t)=\int_0^t u\,e^{-\lambda u}\,du\quad(\lambda>0).
\]
$F(0)=0$, and \eqref{2Equation4} gives
\[
 F'(t)=\frac{pq}{p-q}\bigl(e^{-qt}-e^{-pt}\bigr)
      =\frac{pq}{p-q}\int_q^p t\,e^{-\lambda t}\,d\lambda,
\]
and integration gives
\begin{equation}\label{2Equation6}
 F(t)=\frac{pq}{p-q}\int_q^p F_\lambda(t)\,d\lambda.
\end{equation}
Thus $F\ge0$, $F$ increases on $[0,\infty)$, and $F''>0$ on $(-\infty,0]$.
For $t>0$, $F_\lambda(-t)=\int_0^t u\,e^{\lambda u}\,du>F_\lambda(t)$.

\smallskip
\noindent\textit{Step 1. A lower bound for $6+\mathcal E_{p,q}(X)$.}
Translate the point in the convex hull to the
origin and write $u_i=|x_i|$ and $t_i=\max\{0,-\log u_i\}$ for the remaining three points. If one of
the three angles around the origin exceeded $\pi$, the three points would lie in an open half-plane bounded by a line through the origin, contradicting the choice of that point. Ordering
the gaps as $\phi_1\ge\phi_2\ge\phi_3$, we therefore have
\[
 \phi_1+\phi_2+\phi_3=2\pi,\qquad
 \frac{2\pi}{3}\le\phi_1\le\pi,\qquad
 2\phi_2\ge2\pi-\phi_1\ge\pi.
\]
In the collinear case the gaps are $\pi,\pi,0$. Let $(i_k,j_k)$ denote the pair that subtends
the angle $\phi_k$, for $k=1,2$, and set $s_k=(t_{i_k}+t_{j_k})/2$ and $s=\max\{s_1,s_2\}$. Because $u_{i_k}u_{j_k}\ge
e^{-2s_k}\ge e^{-2s}$, the law of cosines yields
\[
 r_{i_kj_k}^2=(u_{i_k}-u_{j_k})^2+2u_{i_k}u_{j_k}(1-\cos\phi_k)
 \ge2e^{-2s}(1-\cos\phi_k),\quad k=1,2,
\]
and hence $r_{i_1j_1}\ge\sqrt3\,e^{-s}$ and $r_{i_2j_2}\ge\sqrt2\,e^{-s}$. For the pair with $s_k=s$, convexity gives
\[
 F(\log u_{i_k})+F(\log u_{j_k})\ge F(-t_{i_k})+F(-t_{j_k})\ge2F(-s).
\]
Each of the two outer pairs contributes $F(\log r_{i_kj_k})\ge F\bigl((\log r_{i_kj_k})_+\bigr)$, where $t_+=\max\{t,0\}$,
The four pairs are distinct, and $F\ge0$ bounds the remaining pair terms.
Since $F$ increases on $[0,\infty)$,
\begin{equation}\label{2Equation7}
 6+\mathcal E_{p,q}(X)\ge2F(-s)+F\bigl((\log\sqrt2-s)_+\bigr)
                              +F\bigl((\log\sqrt3-s)_+\bigr).
\end{equation}

\smallskip
\noindent\textit{Step 2. The comparison with a square.}
Let $Y$ be a unit square. Its four sides contribute
$4F(0)=0$, while the two diagonals give $6+\mathcal E_{p,q}(Y)=2F(\log\sqrt2)$. If $s=0$, the right-hand side of \eqref{2Equation7} reads $F(\log\sqrt2)+F(\log\sqrt3)>2F(\log\sqrt2)$. If $s\ge\log\sqrt2$, it is at least
$2F(-s)>2F(s)\ge2F(\log\sqrt2)$. In either case, \eqref{2Equation7} gives $\mathcal E_{p,q}(X)>\mathcal E_{p,q}(Y)$.

For $0<s<\log\sqrt2$, let $Y$ be the square of side $e^{-s/2}$. Then
\[
 6+\mathcal E_{p,q}(Y)=4F(-s/2)+2F(\log\sqrt2-s/2).
\]
For every $\lambda>0$, $4F_\lambda(-s/2)=\int_0^s u\,e^{\lambda u/2}\,du<F_\lambda(-s)$, whence
\begin{equation}\label{2Equation8}
 2F_\lambda(-s)-4F_\lambda(-s/2)>F_\lambda(-s)>\frac{s^2}{2}.
\end{equation}
We next estimate the difference between the diagonal terms.

First suppose $0<s\le\log\sqrt{3/2}$. Since $F_\lambda''(t)=(1-\lambda t)e^{-\lambda t}\ge
-e^{-2}$ for $t\ge0$,
\[
\begin{aligned}
 &2F_\lambda(\log\sqrt2-s/2)-F_\lambda(\log\sqrt2-s)-F_\lambda(\log\sqrt3-s)\\
 &\qquad\le\int_{\log\sqrt2-s}^{\log\sqrt2-s/2}
   \bigl[F_\lambda'(t)-F_\lambda'(t+s/2)\bigr]\,dt
 \le\frac{e^{-2}s^2}{4}<\frac{s^2}{2}.
\end{aligned}
\]
For $\log\sqrt{3/2}<s<\log\sqrt2$, we have
$\log\sqrt3-s>\log\sqrt2-s/2$ and $s>\tfrac47\log\sqrt2$. Using $F_\lambda'(t)\le t$ for $t\ge0$,
\[
\begin{aligned}
 &2F_\lambda(\log\sqrt2-s/2)-F_\lambda(\log\sqrt2-s)-F_\lambda(\log\sqrt3-s)\\
 &\qquad\le F_\lambda(\log\sqrt2-s/2)-F_\lambda(\log\sqrt2-s)
 \le\int_{\log\sqrt2-s}^{\log\sqrt2-s/2}t\,dt
 =\frac{s\log\sqrt2-\frac34s^2}{2}<\frac{s^2}{2}.
\end{aligned}
\]
For $0<s<\log\sqrt2$, the difference between the diagonal terms is less than $s^2/2$,
while the side terms increase by more than $s^2/2$ by \eqref{2Equation8}.
Integrating in $\lambda$ and using \eqref{2Equation6} and \eqref{2Equation7} gives
\[
\begin{aligned}
 6+\mathcal E_{p,q}(X)&\ge2F(-s)+F(\log\sqrt2-s)+F(\log\sqrt3-s)\\
 &>4F(-s/2)+2F(\log\sqrt2-s/2)=6+\mathcal E_{p,q}(Y),
\end{aligned}
\]
which proves the comparison with the square of side $e^{-s/2}$.
\end{proof}

\begin{lemma}[Rhombic symmetry]\label{2Lemma3}
Let $p>q>0$. Every strictly convex critical configuration
of $\mathcal E_{p,q}$ is a rhombus. In particular, every global minimizer is a rhombus.
\end{lemma}
\begin{proof}
Let $X$ be a strictly convex critical configuration, and label its vertices in cyclic order.

\smallskip
\noindent\textit{Step 1. The signs of the forces.}
Set $w_{ij}=W'(r_{ij})/r_{ij}$ and $\Omega_{ij}=-w_{ij}$ for $i\ne j$, while the
diagonal entries are defined by
\[
 \Omega_{ii}=\sum_{\substack{1\le j\le4\\j\ne i}}w_{ij},\qquad 1\le i\le4.
\]
Write $x_j=(x_{j,1},x_{j,2})$. The definition of $\Omega_{ii}$ and the criticality of
$X$ for \eqref{2Equation2} give, for $1\le i\le4$,
\[
 \sum_{j=1}^4\Omega_{ij}=\Omega_{ii}-\sum_{\substack{1\le j\le4\\j\ne i}}w_{ij}=0,\qquad
 \sum_{j=1}^4\Omega_{ij}x_{j,k}=\sum_{\substack{1\le j\le4\\j\ne i}}w_{ij}(x_{i,k}-x_{j,k})=0,\quad k=1,2.
\]
Since the four points are not collinear, the vectors $(1,1,1,1)$, $(x_{1,1},\ldots,x_{4,1})$ and
$(x_{1,2},\ldots,x_{4,2})$ are linearly independent. Every row of $\Omega$ lies in their one-dimensional orthogonal complement.

If all off-diagonal entries of $\Omega$ vanish, \eqref{2Equation4} gives $r_{ij}=1$ for every $i<j$.
After translating $x_4$ to the origin, the Gram determinant of $x_1,x_2,x_3$ would satisfy
\[
 \det\begin{pmatrix}
 1&1/2&1/2\\
 1/2&1&1/2\\
 1/2&1/2&1
 \end{pmatrix}=\frac12,
\]
contradicting their linear dependence in $\R^2$. Hence some off-diagonal entry is nonzero.

Let the diagonals meet at
\[
 o=\lambda x_1+(1-\lambda)x_3=\mu x_2+(1-\mu)x_4,\qquad0<\lambda,\mu<1,
\]
and set $(z_1,z_2,z_3,z_4)=(\lambda,-\mu,1-\lambda,\mu-1)$. Then
\[
 \sum_{j=1}^4z_j=0,\qquad\sum_{j=1}^4z_jx_{j,k}=0\qquad(k=1,2).
\]
The orthogonal complement is spanned by $z$, so $\Omega_{ij}=c_i z_j$ for some $c_i\in\R$.
By symmetry, $c_i z_j=c_j z_i$. Since every $z_i$ is nonzero,
the ratios $c_i/z_i$ have a common value $\kappa\ne0$. Thus
\begin{equation}\label{2Equation9}
 \Omega_{ij}=\kappa z_i z_j,\qquad1\le i,j\le4,\qquad z_1,z_3>0,\quad z_2,z_4<0.
\end{equation}
If $\kappa>0$, then \eqref{2Equation9} and \eqref{2Equation4} give $r_{13},r_{24}<1$ together with $r_{ij}>1$ for adjacent pairs.
Choose on each diagonal the endpoint closest to $o$, say $x_i$ and $x_j$. These vertices are adjacent, and
\[
 1<r_{ij}\le|x_i-o|+|x_j-o|\le\frac{r_{13}+r_{24}}2<1,
\]
a contradiction. Thus $\kappa<0$, and \eqref{2Equation9} and \eqref{2Equation4} give
\begin{equation}\label{2Equation10}
 r_{12},r_{23},r_{34},r_{14}<1<r_{13},r_{24}.
\end{equation}

\smallskip
\noindent\textit{Step 2. Perpendicular diagonals, one of them bisected.}
For $0<r<1$ set
\begin{equation}\label{2Equation11}
 f(r)=-\frac{W'(r)}r=\frac{pq}{p-q}r^{-p-2}\bigl(1-r^{p-q}\bigr),\qquad0<r<1.
\end{equation}
Differentiating \eqref{2Equation11} shows that $f$ and $r^2f$ are strictly decreasing on this interval.

Place the intersection of the diagonals at the origin and write
\[
 x_1=ae,\quad x_3=-de,\quad x_2=\alpha e',\quad x_4=-\beta e',\quad
 |e|=|e'|=1,\quad a,d,\alpha,\beta>0,
\]
so that
\[
 (z_1,z_2,z_3,z_4)=\Bigl(\frac{d}{a+d},-\frac{\beta}{\alpha+\beta},\frac{a}{a+d},-\frac{\alpha}{\alpha+\beta}\Bigr).
\]
Set $K=-\kappa/((a+d)(\alpha+\beta))>0$. The side lengths lie in $(0,1)$ by
\eqref{2Equation10}, so \eqref{2Equation9} and \eqref{2Equation11} give
\begin{equation}\label{2Equation12}
 f(r_{12})=Kd\beta,\quad f(r_{23})=Ka\beta,\quad
 f(r_{34})=Ka\alpha,\quad f(r_{14})=Kd\alpha.
\end{equation}
Assume the intersection bisects neither diagonal, and relabel so that $a>d$ and $\alpha>\beta$.
Equation~\eqref{2Equation12}, together with the monotonicity of $f$ and of $r^2f$, gives
\[
 \begin{gathered}
 r_{12}>r_{23},\qquad r_{14}>r_{34},\qquad dr_{12}^2<ar_{23}^2,\qquad dr_{14}^2<ar_{34}^2,\\
 r_{12}>r_{14},\qquad r_{23}>r_{34},\qquad \beta r_{12}^2<\alpha r_{14}^2,\qquad \beta r_{23}^2<\alpha r_{34}^2.
 \end{gathered}
\]
For example, the law of cosines gives
\[
 dr_{12}^2-ar_{23}^2=(a-d)(ad-\alpha^2)-4ad\alpha\,e\cdot e',
\]
and the mixed terms cancel in the sums
\[
 \begin{aligned}
 0&>\beta\bigl(dr_{12}^2-ar_{23}^2\bigr)+\alpha\bigl(dr_{14}^2-ar_{34}^2\bigr)
 =(a-d)(\alpha+\beta)(ad-\alpha\beta),\\
 0&>d\bigl(\beta r_{12}^2-\alpha r_{14}^2\bigr)+a\bigl(\beta r_{23}^2-\alpha r_{34}^2\bigr)
 =(\alpha-\beta)(a+d)(\alpha\beta-ad).
 \end{aligned}
\]
Since $a>d$ and $\alpha>\beta$, the first identity gives $ad<\alpha\beta$ and the second gives $\alpha\beta<ad$, a contradiction.
The intersection bisects at least one diagonal. Relabel so that $a=d$.
Equation~\eqref{2Equation12} and the strict decrease of $f$ give $r_{12}=r_{23}$.
Since $r_{23}^2-r_{12}^2=4a\alpha\,e\cdot e'$, we obtain $e\cdot e'=0$.

\smallskip
\noindent\textit{Step 3. The second diagonal is bisected as well.}
Suppose $\alpha\ne\beta$. Interchanging $x_2$ and $x_4$ if necessary, assume $\alpha>\beta$. Write $c=a^2$, $u=c+\alpha^2=r_{12}^2$ and $v=c+\beta^2=r_{14}^2$. By \eqref{2Equation10}
we have $0<v<u<1$, and $\alpha+\beta>1$. Then
\[
 \beta u(1-v)-\alpha v(1-u)=(\alpha-\beta)(uv+\alpha\beta-c),
\]
with
\[
 uv+\alpha\beta-c=(c-\alpha\beta)^2+c\bigl[(\alpha+\beta)^2-1\bigr]+\alpha\beta>0,
\]
so that $\alpha(1-u)/u<\beta(1-v)/v$. For $s>0$, both $t^{-s/2}$ and
\[
 \frac{-\log t}{1-t}=\int_0^1\frac{d\tau}{1-\tau+\tau t}
\]
strictly decrease on $(0,1)$. Multiplying the preceding inequality by these positive factors and using $u>v$ gives
\[
 \alpha(-\log u)\,u^{-(s+2)/2}<\beta(-\log v)\,v^{-(s+2)/2}.
\]
For $0<r<1$, \eqref{2Equation11} also implies
\[
 f(r)=\frac{pq}{p-q}\int_q^p(-\log r)\,r^{-s-2}\,ds.
\]
Substituting $u=r_{12}^2$ and $v=r_{14}^2$ into the previous inequality and integrating over $s\in(q,p)$
yields $\alpha f(r_{12})<\beta f(r_{14})$, which contradicts $\alpha f(r_{12})=Kd\alpha\beta=\beta f(r_{14})$ from \eqref{2Equation12}. Hence
$\alpha=\beta$.

Together with Step~2, this gives $a=d$ and $e\cdot e'=0$, and therefore
\[
 r_{12}^2=r_{23}^2=r_{34}^2=r_{14}^2=a^2+\alpha^2,
\]
so $X$ is a rhombus.

Every global minimizer is strictly convex by Lemma~\ref{2Lemma2} and is critical because
$\mathcal E_{p,q}$ is smooth on the open set of distinct configurations. The preceding argument therefore applies.
\end{proof}

\section{The square and rhombic phases}\label{Section3}

By Lemma~\ref{2Lemma3}, every minimizer is a rhombus for $b>0$.
We prove uniqueness, locate the square region $D_{p,q}\ge1$, and study the
angle as a function of $(p,q)$, including the four-point Heitmann--Radin limit.

A rhombus of side length $a>0$ and angle $0<\theta\le\pi/2$ has diagonals
$2a\sin(\theta/2)$ and $2a\cos(\theta/2)$. For $s>0$, set
\begin{equation}\label{3Equation1}
 S_s(\theta)=4+\bigl(2\sin(\theta/2)\bigr)^{-s}
                  +\bigl(2\cos(\theta/2)\bigr)^{-s}.
\end{equation}
\begin{proposition}[Minimization among rhombi]\label{3Proposition1}
For $p>q>0$ and $b>0$, the minimizing rhombus is unique up to congruence
and relabelling, and its side length is uniquely determined. It is a square if
$D_{p,q}\ge1$, where $D_{p,q}$ is defined in \eqref{1Equation4}.
If $D_{p,q}<1$, its acute angle is
$\theta_{p,q}=2\arctan(e^{-L_{p,q}})\in(\pi/3,\pi/2)$,
where $L_{p,q}$ is the unique solution in $(0,\tfrac12\log3)$ of
\begin{equation}\label{3Equation2}
 \frac{1-e^{-(p+2)L}}{1-e^{-(q+2)L}}
 =\frac{1+\left(\frac4{1+e^{2L}}\right)^{(p+2)/2}}
        {1+\left(\frac4{1+e^{2L}}\right)^{(q+2)/2}}.
\end{equation}
\end{proposition}

To prove Proposition~\ref{3Proposition1}, fix $\theta$ and write
$S_s=S_s(\theta)$. By \eqref{1Equation3} and \eqref{3Equation1}, the rhombus has energy
$S_pa^{-p}-bS_qa^{-q}$. Differentiating with respect to $a$ gives
\[
 a\frac{d}{da}\left(S_pa^{-p}-bS_qa^{-q}\right)
 =a^{-p}\left(bqS_qa^{p-q}-pS_p\right).
\]
The factor in parentheses increases strictly from $-pS_p$ to $+\infty$. Hence
\begin{equation}\label{3Equation3}
 \min_{a>0}\bigl\{S_pa^{-p}-bS_qa^{-q}\bigr\}
 =-\varepsilon_b\frac{S_q^{p/(p-q)}}{S_p^{q/(p-q)}},
 \qquad a^{p-q}=\frac{pS_p}{bqS_q}.
\end{equation}
Thus the minimizing angle maximizes $S_q^p/S_p^q$.
By \eqref{2Equation3} and Lemmas~\ref{2Lemma1} and~\ref{2Lemma3}, a rhombus attains
the global minimum. Scaling this minimizer by $\ell_b^{-1}$ and applying
\eqref{2Equation10} gives
\[
 1<\frac{2a\sin(\theta/2)}{\ell_b}<2\sin(\theta/2)
 \quad\Longrightarrow\quad \theta>\pi/3.
\]
Set $L=\log\cot(\theta/2)=\operatorname{arctanh}(\cos\theta)$. Then
\[
 S_s(\theta)=4+2\left(\frac{\cosh L}{2}\right)^{s/2}\cosh(sL/2),
 \qquad
 \frac{dS_s}{dL}
 =\frac{s}{\cosh L}\left(\frac{\cosh L}{2}\right)^{s/2}
       \sinh\bigl((s+2)L/2\bigr).
\]
Define $R(L)=R_1(L)+R_2(L)$ for $L>0$, where
\begin{equation}\label{3Equation4}
 \begin{aligned}
 R_1(L)&=\left(\frac{\cosh L}{2}\right)^{(p-q)/2}
           \frac{\sinh\bigl((p+2)L/2\bigr)}{\sinh\bigl((q+2)L/2\bigr)},\\
 R_2(L)&=\left(\frac{\cosh L}{2}\right)^{(p+2)/2}
           \frac{\sinh\bigl((p-q)L/2\bigr)}{\sinh\bigl((q+2)L/2\bigr)}.
 \end{aligned}
\end{equation}
The addition formulas for $\sinh$ give
\[
 \sinh\bigl((q+2)L/2\bigr)\cosh(pL/2)
   -\sinh\bigl((p+2)L/2\bigr)\cosh(qL/2)
 =-\cosh L\,\sinh\bigl((p-q)L/2\bigr),
\]
and therefore
\begin{equation}\label{PL1}
 \frac{d}{dL}\log\frac{S_q^p}{S_p^q}
 =\frac{4pq}{\cosh L\,S_pS_q}\left(\frac{\cosh L}{2}\right)^{q/2}
             \sinh\frac{(q+2)L}{2}\,[1-R(L)].
\end{equation}

\begin{lemma}[Monotonicity of $R(L)$]\label{3Lemma2}
The function $R$ is strictly increasing on $(0,\tfrac12\log3)$. Its
endpoint values satisfy
\[
 R(0)=D_{p,q},\qquad R(\tfrac12\log3)>1,
\]
where $R(0)$ is defined by continuity.
\end{lemma}
\begin{proof}
Set $z=sL>0$. Since $\sinh z\cosh z>z$ and $z\coth z>1$,
\[
 0<\frac{\sinh z\cosh z-z}{\sinh^2z}
 =\frac{\partial}{\partial s}\bigl[s\coth(sL)\bigr]
 =z-\coth z\bigl(z\coth z-1\bigr)<z.
\]
Taking logarithmic derivatives in \eqref{3Equation4} gives
\[
 \begin{aligned}
 \frac{R_1'}{R_1}
 &=\frac{p-q}{2}\tanh L
   +\int_{(q+2)/2}^{(p+2)/2}\frac{\partial}{\partial s}
                                      \bigl[s\coth(sL)\bigr]\,ds
 >\frac{p-q}{2}\tanh L,\\
 \frac{R_2'}{R_2}
 &=\frac{p+2}{2}\tanh L
   +\int_{(q+2)/2}^{(p-q)/2}\frac{\partial}{\partial s}
                                      \bigl[s\coth(sL)\bigr]\,ds.
 \end{aligned}
\]
If $p-q\ge q+2$, both logarithmic derivatives are positive, so $R'=R_1'+R_2'>0$.

If $p-q<q+2$,
\[
 \begin{aligned}
 \frac{R_2'}{R_2}
 &>-\int_{(p-q)/2}^{(q+2)/2}sL\,ds
   =-\frac{(2q+2-p)(p+2)L}{8},\\
 \frac{R_2}{R_1}
 &\le\left(\frac{\cosh L}{2}\right)^{(q+2)/2}\frac{p-q}{p+2},
 \end{aligned}
\]
where the second inequality follows because $\sinh z/z=\int_0^1\cosh(tz)\,dt$ increases on $(0,\infty)$. For
$0<L<\tfrac12\log3$,
\[
 \frac12\cosh L<\frac1{\sqrt3},\qquad
 \tanh L=\int_0^L\operatorname{sech}^2t\,dt\ge\frac{3L}{4}.
\]
Hence, whenever $p-q<q+2$,
\[
 \frac{2q+2-p}{2}\left(\frac{\cosh L}{2}\right)^{(q+2)/2}
 <\frac{q+2}{2}\,3^{-(q+2)/4}
 \le \sup_{x>0}x\,3^{-x/2}
 =\frac{2}{e\log3}<1,
\]
and combining these bounds for $R_1'/R_1$, $R_2'/R_2$ and $R_2/R_1$ gives
\begin{equation}\label{3Equation5}
 \begin{aligned}
 \frac{R'}{R_1}
 &=\frac{R_1'}{R_1}+\frac{R_2}{R_1}\frac{R_2'}{R_2}\\
 &>\frac{p-q}{2}\left[
    \tanh L-\frac{(2q+2-p)L}{4}
                   \left(\frac{\cosh L}{2}\right)^{(q+2)/2}\right]
 >\frac{(p-q)L}{8}.
 \end{aligned}
\end{equation}
Since $R_1(L)>0$, this proves $R'(L)>0$ in the second case.

Letting $L\downarrow0$ in \eqref{3Equation4} and using $\sinh t/t\to1$, we obtain
\[
 R(0)=2^{-(p-q)/2}\left[\frac{p+2}{q+2}
      +2^{-(q+2)/2}\frac{p-q}{q+2}\right]=D_{p,q},
\]
where the last equality follows from \eqref{1Equation4}. For the other endpoint, set
$x=\cos\theta=\tanh L$. Then
\[
 S_s(\arccos x)=4+[2(1-x)]^{-s/2}+[2(1+x)]^{-s/2}.
\]
At $L=\tfrac12\log3$, we have $x=1/2$, so
\[
 \left.\frac1s\frac{d}{dx}\log S_s(\arccos x)\right|_{x=1/2}
 =\frac{1-3^{-s/2-1}}{5+3^{-s/2}},
\]
and
\[
 \frac{d}{ds}\frac{1-3^{-s/2-1}}{5+3^{-s/2}}
 =\frac{4\log3}{3}\frac{3^{-s/2}}{(5+3^{-s/2})^2}>0.
\]
At $L=\tfrac12\log3$, since $p>q$,
\[
 \frac{d}{dL}\log\frac{S_q^p}{S_p^q}
 =pq(1-x^2)
 \left[\frac{1-3^{-q/2-1}}{5+3^{-q/2}}
      -\frac{1-3^{-p/2-1}}{5+3^{-p/2}}\right]<0.
\]
The prefactor in \eqref{PL1} is positive for $L>0$, so the derivative has
the sign of $1-R(L)$. The endpoint inequality therefore gives
$R(\tfrac12\log3)>1$.
\end{proof}

\begin{proof}[Proof of Proposition~\ref{3Proposition1}]
By \eqref{PL1}, the derivative of $\log(S_q^p/S_p^q)$ has the sign of
$1-R(L)$. Lemma~\ref{3Lemma2} therefore gives a unique maximizing value of $L$,
\[
 L=\begin{cases}
 0, & D_{p,q}\ge1,\\
 L_{p,q}\in(0,\tfrac12\log3),\quad R(L_{p,q})=1, & D_{p,q}<1.
 \end{cases}
\]
For $s,L>0$, set
\[
 K_s(L)=\log(1-e^{-sL})
         -\log\left(1+\left(\frac4{1+e^{2L}}\right)^{s/2}\right).
\]
The identity $1-e^{-sL}=2e^{-sL/2}\sinh(sL/2)$ gives
\begin{equation}\label{3Equation6}
 e^{K_{p+2}(L)-K_{q+2}(L)}-1
 =\frac{R(L)-1}{1+\left(\frac{1+e^{2L}}4\right)^{(p+2)/2}}.
\end{equation}
Thus $R(L)=1$ is equivalent to \eqref{3Equation2}. Differentiating
\eqref{3Equation6} at the root and using $R'(L)>0$ from Lemma~\ref{3Lemma2} gives
\[
 \partial_L\bigl(K_{p+2}-K_{q+2}\bigr)
 =\frac{R'(L)}{1+\left(\frac{1+e^{2L}}4\right)^{(p+2)/2}}>0.
\]

The map $L\mapsto2\arctan(e^{-L})$ is one-to-one, so the maximizing
value above determines a unique minimizing angle, equal to $\pi/2$ exactly
when $D_{p,q}\ge1$. Equation~\eqref{3Equation3} gives a unique minimizing
side length. The rhombus is therefore unique up to congruence and relabelling.
\end{proof}

Lemmas~\ref{2Lemma1} and~\ref{2Lemma3} and Proposition~\ref{3Proposition1} prove
parts \textup{(ii)--(iii)} of Theorem~\ref{1Theorem1}. Equation~\eqref{3Equation3}
gives $a_b$ and $e_{p,q,b}$ at $\theta=\theta_{p,q}$.

To locate the curve $D_{p,q}=1$, define
\begin{equation}\label{3Equation7}
 H(s)=\log s-\log\bigl(1+2^{s/2}\bigr),\qquad s>0.
\end{equation}
\begin{lemma}[The square region]\label{3Lemma3}
There is a unique $q_0\in(1,2)$ with $H'(q_0+2)=0$. For
$0<q<q_0$, the equation $D_{p,q}=1$ has a unique solution $p=p_{\rm c}(q)>q$. The function $p_{\rm c}$ is real
analytic and strictly decreasing, and its limits are
\[
 \lim_{q\downarrow0}p_{\rm c}(q)=4,\qquad
 \lim_{q\uparrow q_0}p_{\rm c}(q)=q_0.
\]
For $p>q>0$ and $b>0$, the square region in Figure~\ref{1Figure2} is
\[
 \theta_{p,q}=\pi/2
 \quad\Longleftrightarrow\quad
 0<q<q_0 \text{ and }q<p\le p_{\rm c}(q).
\]
\end{lemma}
\begin{proof}
Equation \eqref{1Equation4} gives
\[
 (p+2)\bigl(1+2^{(q+2)/2}\bigr)-(q+2)\bigl(1+2^{(p+2)/2}\bigr)
 =(q+2)2^{(p+2)/2}\bigl(D_{p,q}-1\bigr),
\]
so \eqref{3Equation7} gives
$\operatorname{sgn}\bigl(H(p+2)-H(q+2)\bigr)=\operatorname{sgn}(D_{p,q}-1)$. Moreover,
\[
 H'(s)=\frac1s-\frac{\log2}{2}\frac{2^{s/2}}{1+2^{s/2}},
 \qquad
 H''(s)=-\frac1{s^2}
 -\frac{(\log2)^2}{4}\frac{2^{s/2}}{(1+2^{s/2})^2}<0.
\]
Since $H'(3)>0>H'(4)$ and $H''<0$, there is a unique solution
$s=q_0+2\in(3,4)$ of $H'(s)=0$. The function $H$ increases strictly up to
$q_0+2$ and then decreases strictly to $-\infty$. For $0<q<q_0$, we have
$H(q+2)<H(q_0+2)$. The decreasing branch therefore contains exactly one point
$p_{\rm c}(q)+2$ with $H(p_{\rm c}(q)+2)=H(q+2)$, and the increasing branch gives
no solution with $q<p\le q_0$. The sign identity gives a unique root of
$D_{p,q}=1$ with $p>q$.

Since $H$ is real analytic and $H'(p_{\rm c}(q)+2)<0$, the implicit function
theorem shows that $p_{\rm c}$ is real analytic and
\[
 p_{\rm c}'(q)=\frac{H'(q+2)}{H'(p_{\rm c}(q)+2)}<0.
\]
As $q\downarrow0$, we have $H(q+2)\to H(2)=H(6)$. Continuity of the inverse
of the decreasing branch gives $p_{\rm c}(q)+2\to6$. As $q\uparrow q_0$, we have
$H(q+2)\to H(q_0+2)$, and the same inverse gives $p_{\rm c}(q)+2\to q_0+2$.

For $q\ge q_0$ and $p>q$, strict decrease of $H$ gives $H(p+2)<H(q+2)$. For
$0<q<q_0$, strict increase on $[q+2,q_0+2]$ and strict decrease on
$[q_0+2,\infty)$ give $H(p+2)\ge H(q+2)$ exactly when $q<p\le p_{\rm c}(q)$.
Proposition~\ref{3Proposition1} and the sign identity therefore give
\[
 \theta_{p,q}=\pi/2
 \iff D_{p,q}\ge1
 \iff H(p+2)\ge H(q+2)
 \iff 0<q<q_0\text{ and }q<p\le p_{\rm c}(q).
\]
\end{proof}

By Lemma~\ref{3Lemma3}, the strict square region above $q=1$ is
\[
 \{(p,q):p>q>1,\ D_{p,q}>1\}
 =\{(p,q):1<q<q_0,\ q<p<p_{\rm c}(q)\}.
\]
Equation~\eqref{1Equation4} gives $D_{2,4/3}>1>D_{5/2,1}$. Hence
Lemma~\ref{3Lemma3} gives $p_{\rm c}(1)<5/2$. Since $p_{\rm c}$ decreases,
$D_{p,q}<1$ for $p>q>1$ with $p\ge5/2$.

\begin{proof}[Proof of part \textup{(iv)} of Theorem~\ref{1Theorem1}]
\textup{(a)} Fix $0<L<\tfrac12\log3$ and set $t=\tfrac12\log\bigl(4/(1+e^{2L})\bigr)$.
Then
\[
 \partial_s^2 K_s(L)
 =-\frac{L^2e^{sL}}{(e^{sL}-1)^2}
  -\frac{t^2e^{st}}{(1+e^{st})^2}<0.
\]
By Proposition~\ref{3Proposition1} and \eqref{3Equation6}, $K_{p+2}(L_{p,q})=K_{q+2}(L_{p,q})$. Strict concavity in $s$ then implies
\[
 \left.\partial_sK_s\right|_{s=q+2}
 >\frac{K_{p+2}-K_{q+2}}{p-q}=0
 >\left.\partial_sK_s\right|_{s=p+2}.
\]
The positive $L$-derivative of $K_{p+2}-K_{q+2}$ at the root, proved in
Proposition~\ref{3Proposition1}, allows us to apply the implicit function
theorem. Thus $L_{p,q}$ is real analytic, and differentiating the defining equation gives
\[
 \partial_p L_{p,q}
 =-\frac{\left.\partial_sK_s\right|_{s=p+2}}
         {\partial_L(K_{p+2}-K_{q+2})}>0,\qquad
 \partial_q L_{p,q}
 =\frac{\left.\partial_sK_s\right|_{s=q+2}}
         {\partial_L(K_{p+2}-K_{q+2})}>0.
\]
By Proposition~\ref{3Proposition1}, $\theta_{p,q}=2\arctan(e^{-L_{p,q}})$.
Since $\partial_L\bigl[2\arctan(e^{-L})\bigr]=-\operatorname{sech}L$, both angle
derivatives are negative.

\textup{(b)} Fix $q>0$. By \eqref{1Equation4}, $D_{p,q}\to0$ as $p\to\infty$; choose $p_1>q$ so that $D_{p,q}<1$ for $p\ge p_1$.
Part \textup{(a)} then gives $L_{p,q}\ge L_{p_1,q}>0$ on this range. Equation \eqref{3Equation2} reads
\[
 \left(\frac4{1+e^{2L_{p,q}}}\right)^{(p+2)/2}
 =\frac{1-e^{-(p+2)L_{p,q}}}{1-e^{-(q+2)L_{p,q}}}
 \left[1+\left(\frac4{1+e^{2L_{p,q}}}\right)^{(q+2)/2}\right]-1
 =O(1).
\]
Since $4/(1+e^{2L})$ exceeds $1$ for $L<\tfrac12\log3$, boundedness of the
left-hand side gives $L_{p,q}\to\tfrac12\log3$. Passing to the limit on the
right-hand side gives
\[
 \left(\frac4{1+e^{2L_{p,q}}}\right)^{(p+2)/2}
 \longrightarrow\frac{1+3^{-(q+2)/2}}{1-3^{-(q+2)/2}}.
\]
At $L=\tfrac12\log3$ the derivatives of $2\arctan(e^{-L})$ and $\tfrac12\log\bigl(4/(1+e^{2L})\bigr)$ are $-\frac{\sqrt3}{2}$ and $-\frac34$, hence
\[
 \lim_{L\uparrow\tfrac12\log3}
 \frac{2\arctan(e^{-L})-\pi/3}{\tfrac12\log\bigl(4/(1+e^{2L})\bigr)}
 =\frac2{\sqrt3}.
\]
Multiplying this limit by the logarithm of the preceding limit proves
part \textup{(iv)(b)}.

\textup{(c)} Take $p=2q$ with $q\ge2$; Lemma~\ref{3Lemma3} gives $D_{2q,q}<1$, and part \textup{(a)} gives $L_{2q,q}\ge L_{4,2}>0$.
Write $L=L_{2q,q}$, $L_0=\tfrac12\log3$ and $T=\frac4{1+e^{2L}}>1$. Equation \eqref{3Equation2} becomes
\[
 \frac{T^{q/2}-1}{1+T^{-(q+2)/2}}
 =\frac{e^{-(q+2)L}\bigl(1-e^{-qL}\bigr)}{1-e^{-(q+2)L}}.
\]
Because $\log T=\int_L^{L_0}\frac{2e^{2t}}{1+e^{2t}}\,dt\ge L_0-L$,
\[
 0<q(L_0-L)\le q\log T\le 2\bigl(T^{q/2}-1\bigr)=O\bigl(e^{-(q+2)L}\bigr)=o(1).
\]
Thus $e^{-(q+2)L}=3^{-(q+2)/2}e^{(q+2)(L_0-L)}=O(3^{-q/2})$, so
$q(L_0-L)=O(3^{-q/2})$. Expanding both sides of \eqref{3Equation2} gives
\[
 \begin{aligned}
 \frac{T^{q/2}-1}{1+T^{-(q+2)/2}}
 &=\frac q4\log T+O(3^{-q})
  =\frac{3q}{8}(L_0-L)+O(3^{-q}),\\
 \frac{e^{-(q+2)L}\bigl(1-e^{-qL}\bigr)}{1-e^{-(q+2)L}}
 &=3^{-(q+2)/2}+O(3^{-q}),
 \end{aligned}
\]
Hence $L_0-L=\frac{8}{9q}3^{-q/2}+O(3^{-q}/q)$. The angle formula in
Proposition~\ref{3Proposition1} gives
\[
 \theta_{2q,q}-\frac\pi3
 =\frac{\sqrt3}{2}(L_0-L)+O\bigl((L_0-L)^2\bigr)
 =\frac{4}{3\sqrt3}q^{-1}3^{-q/2}+O(3^{-q}/q).
\]
At $\theta=\theta_{2q,q}$, Taylor expansion gives
\[
 q\log\bigl(2\sin(\theta/2)\bigr)
 =\frac23\,3^{-q/2}+O(3^{-q}),\qquad
 q\log\bigl(2\cos(\theta/2)/\sqrt3\bigr)
 =O(3^{-q/2}),
\]
and substitution into \eqref{3Equation1} gives
\[
 S_q(\theta)=5+\frac13\,3^{-q/2}+O(3^{-q}),\qquad
 S_{2q}(\theta)=5-\frac43\,3^{-q/2}+O(3^{-q}).
\]
For $b=2$, equation \eqref{3Equation3} then reads
\[
 q\log a_2=\log\frac{S_{2q}(\theta)}{S_q(\theta)}
 =-\frac13\,3^{-q/2}+O(3^{-q}),
\]
and in particular $a_2\to1$.

\textup{(d)} Write $R_{p,q}$ for $R$ in \eqref{3Equation4}. Dividing the numerator and denominator of each hyperbolic quotient by $L$
gives a joint analytic extension near $(p^*,q^*,0)$, even in $L$. The expansion
\[
 (\cosh L)^\gamma\frac{\sinh(uL)}{\sinh(vL)}
 =\frac uv\left[1+\left(\frac{\gamma}{2}+\frac{u^2-v^2}{6}\right)
                      L^2+O(L^4)\right]
\]
gives
\[
 \begin{aligned}
 R_{p,q}''(0)
 &=\frac{(p-q)(p+q+10)}{12}\,R_{p,q}(0)
   +\frac{(q+2)(4-2p+q)}{12}\,R_2(0)\\
 &=\frac{p-q}{12}\left[(p+q+10)D_{p,q}
                       +(4-2p+q)2^{-(p+2)/2}\right].
 \end{aligned}
\]
On the curve $D_{p,q}=1$ the bracket equals $p(1-2^{-p/2})+q+10+(q+4)2^{-(p+2)/2}>0$.

Let $(p,q)\to(p^*,q^*)$ through $D_{p,q}<1$. Proposition~\ref{3Proposition1} and
\eqref{3Equation6} give $R_{p,q}(L_{p,q})=1$. If a subsequence satisfied
$L_{p,q}\to L^*>0$, Lemma~\ref{3Lemma2} and continuity would give
$1=R_{p^*,q^*}(L^*)>R_{p^*,q^*}(0)=1$, a contradiction. Thus $L_{p,q}\to0$.
The even Taylor expansion is uniform near $(p^*,q^*)$ and gives
\[
 1-D_{p,q}=\tfrac12 R_{p,q}''(0)L_{p,q}^{2}+O(L_{p,q}^{4}),
 \qquad
 \pi/2-2\arctan(e^{-L_{p,q}})=L_{p,q}+O(L_{p,q}^{3}).
\]
Therefore
\[
 \frac{\pi/2-\theta_{p,q}}{\sqrt{1-D_{p,q}}}
 =\frac{1+O(L_{p,q}^{2})}
 {\sqrt{\tfrac12R_{p,q}''(0)+O(L_{p,q}^{2})}}
 \longrightarrow\sqrt{\frac{2}{R_{p^*,q^*}''(0)}},
\]
and substitution of the formula for $R''(0)$ gives the stated constant.
\end{proof}

\begin{proof}[Proof of Remark~\ref{1Remark2}]
Part \textup{(iv)(c)} of the proof of Theorem~\ref{1Theorem1} gives $a_2\to1$ and
$\theta_{2q,q}\to\pi/3$. Substitution in part \textup{(ii)} gives the stated matrix.
Its pair distances are
\[
 r_{12}=r_{23}=r_{34}=r_{14}=r_{24}=1,\qquad r_{13}=\sqrt3.
\]
The Gram determinant argument in the proof of Lemma~\ref{2Lemma3} rules out
six unit pair distances in the plane. Thus every four-point configuration
with finite Heitmann--Radin energy has at most five unit bonds, and the
displayed limit is a minimizer.
\end{proof}

\section{The classical Lennard-Jones case}\label{Section4}

Theorem~\ref{1Theorem1} gives a sharp inequality for the six pair distances.
For $(p,q)=(12,6)$, we determine the minimizing angle, side length, and energy explicitly.

\begin{corollary}[Sharp scale-free inequality]\label{4Corollary1}
For every configuration $X$ and $p>q>0$,
\begin{equation}\label{4Equation1}
 \frac{\left(\sum_{1\le i<j\le4}|x_i-x_j|^{-q}\right)^p}
      {\left(\sum_{1\le i<j\le4}|x_i-x_j|^{-p}\right)^q}
 \le\frac{S_q(\theta_{p,q})^p}{S_p(\theta_{p,q})^q}.
\end{equation}
Equality holds if and only if $X$ is similar to the minimizing
rhombus in Theorem~\ref{1Theorem1}.
\end{corollary}

\begin{proof}
Fix $b>0$. By Theorem~\ref{1Theorem1}\textup{(ii)} and the scale minimization~\eqref{3Equation3},
\[
 \begin{aligned}
 \min_{t>0}E_{p,q,b}(tX)
 &=-\varepsilon_b
 \left[\frac{\left(\sum_{1\le i<j\le4}|x_i-x_j|^{-q}\right)^p}
            {\left(\sum_{1\le i<j\le4}|x_i-x_j|^{-p}\right)^q}\right]^{1/(p-q)}\\
 &\ge e_{p,q,b}
 =-\varepsilon_b
 \left[\frac{S_q(\theta_{p,q})^p}{S_p(\theta_{p,q})^q}\right]^{1/(p-q)}.
 \end{aligned}
\]
Since $\varepsilon_b>0$, this proves \eqref{4Equation1}.
Lemma~\ref{2Lemma3} and Proposition~\ref{3Proposition1} give
\[
 \min_{t>0}E_{p,q,b}(tX)=e_{p,q,b}
 \quad\Longleftrightarrow\quad
 X\text{ is similar to }X_{p,q,b}^*.
\]
\end{proof}

\begin{proposition}[Classical Lennard--Jones potential]\label{4Proposition2}
Let $b>0$. The minimizing angle is $\theta_{12,6}=\arccos\sqrt{u^*}$,
where $u^*$ is the unique root in $(0,1/4)$ of
\begin{equation}\label{4Equation2}
 397-1585u-336u^2+1008u^3-512u^4=0.
\end{equation}
The minimizing side length and energy satisfy
\begin{equation}\label{4Equation3}
 \begin{aligned}
 a_b^6&=\frac{128(1-u^*)^6+1+15u^*+15u^{*2}+u^{*3}}
              {4b(1-u^*)^3[16(1-u^*)^3+1+3u^*]},\\
 e_{12,6,b}&=-\frac{b^2[16(1-u^*)^3+1+3u^*]^2}
                    {2[128(1-u^*)^6+1+15u^*+15u^{*2}+u^{*3}]}.
 \end{aligned}
\end{equation}
\end{proposition}

\begin{proof}
By \eqref{1Equation4}, $D_{12,6}=115/512<1$. Theorem~\ref{1Theorem1}\textup{(iii)} therefore gives a minimizing
angle in $(\pi/3,\pi/2)$. Set $u=\cos^2\theta\in(0,1/4)$. Then
\[
 4\sin^2(\theta/2)=2(1-\sqrt u),\qquad
 4\cos^2(\theta/2)=2(1+\sqrt u),
\]
and expansion of \eqref{3Equation1} gives
\[
 S_6=\frac{16(1-u)^3+1+3u}{4(1-u)^3},\qquad
 S_{12}=\frac{128(1-u)^6+1+15u+15u^2+u^3}{32(1-u)^6}.
\]
Differentiating, we obtain
\[
 \frac{d}{du}\log\frac{S_6^2}{S_{12}}
 =\frac{3(1-u)^2(397-1585u-336u^2+1008u^3-512u^4)}
        {[16(1-u)^3+1+3u][128(1-u)^6+1+15u+15u^2+u^3]}.
\]
By \eqref{3Equation3}, the minimizing angle maximizes $S_6^2/S_{12}$.
Thus $u=\cos^2\theta_{12,6}$ satisfies \eqref{4Equation2}.
On $[0,1/4]$, the derivative of the polynomial is at most $-1585+3024/16<0$, so this root is unique.
By \eqref{3Equation3},
\[
 a_b^6=\frac{2S_{12}}{bS_6},\qquad
 e_{12,6,b}=-\frac{b^2S_6^2}{4S_{12}},
\]
and substitution of $S_6$ and $S_{12}$ at $u=u^*$ gives \eqref{4Equation3}.
\end{proof}

For $b=2$,
\[
 u^*\approx0.2459296640551374,\qquad a_2\approx0.9980123196944918,\qquad
 e_{12,6,2}\approx-5.073420858462792.
\]

\bigskip
\noindent
{\bf Acknowledgements.}
 The research of S. Luo is partially supported by the National Natural Science Foundation of China (NSFC) under Grant Nos. 12261045 and 12001253, and by the Jiangxi Jieqing Fund under Grant No. 20242BAB23001.

\vskip0.1in
{\bf Statements and Declarations: there is no conflict of interest.}

{\bf Data availability: the manuscript has no associated data.}

\end{document}